# Modeling Infectious Diseases: From SIR Models to Diffusion-Based Approaches and Numerical Solutions

Ayesha Baig[a], Zhouxin li[a]

[a] *Department of mathematics and Statistics, Central South University, Changsha 410083, PR China*

**Abstract**

As global living standards improve and medical technology advances, many infectious diseases have been effectively controlled. However, certain diseases, such as the recent COVID-19 pandemic, continue to pose significant threats to public health. This paper explores the evolution of infectious disease modeling, from early ordinary differential equation-based models like the SIR framework to more complex reaction-diffusion models that incorporate both temporal and spatial dynamics. The study highlights the importance of numerical methods, such as the Runge-Kutta method, implicit-explicit time-discretization techniques, and finite difference methods, in solving these models. By analyzing the development and application of these methods, this research underscores their critical role in predicting disease spread, informing public health strategies, and mitigating the impact of future pandemics.

## 1. Introduction

The study of infectious disease spread began in the early 20th century. In 1906, Hamer [1] developed a time-discrete model for measles spread and analyzed it, offering effective control suggestions. In 1911, Ross [2] established a differential equation model to analyze malaria transmission dynamics, providing a theoretical foundation for controlling its spread. In 1926, Kermack and McKendrick [3] introduced the SIR model to analyze the spread of plague in London and Bombay. Since the mid-20th century, the development of mathematical models for infectious diseases has accelerated. In 1957, Bailey [4] published *Mathematical Epidemiology*, which greatly advanced the study of disease spread. In 1980, Aronsson and Mellander [5] incorporated seasonal variations into a gonorrhea reaction-diffusion model and analyzed its stability and threshold conditions based on Perron-Frobenius eigenvalues. In 2008, Liu and Takeuchi [6] proposed an optimal control model for viral infections based on the SVIR model, describing its dynamic properties.

Since the outbreak of the novel coronavirus in 2020, there has been further exploration into infectious disease models. Tang [7] and others developed an improved SEIR model with 8 groups to study epidemic dynamics; Bogoch [8] examined how travel factors influence the spread; Zhang [9] applied statistical methods to analyze virus diffusion; Giordano [10,11,12] considered the impact of asymptomatic carriers on virus spread, establishing dynamic models for epidemic propagation.

In 1895, Runge [13] introduced a new numerical method for solving ordinary differential equations in his classic work *Numerical Solutions of Ordinary Differential Equations*, which launched the study of the Runge-Kutta method. In 1901, Kutta, building on Runge and Heun's work [13,14], presented the classical Runge-Kutta method formulation. Nystrom [15] proposed a 6th-order 5-stage explicit method in 1925, and Huta [16,17] introduced an 8th-order 6-stage method in 1956, improving the numerical simulation accuracy. Later, in 1972, Verner [18] developed an 11th-order 8-stage explicit method. As computing power increased, the Runge-Kutta method became more widely applied, offering greater stability and accuracy.

In 1995, Ascher [19] proposed implicit-explicit linear multistep methods and derived conditions for stability and error estimates. Over time, this approach has been applied to stiff differential equations. Frank [20] studied the stability of these methods in 1997, and other scholars have applied them to solve time-dependent partial differential equations. In 2011, Li [21] introduced an implicit-explicit predictor-corrector algorithm and analyzed its stability. In 2013, Akrivis [22] applied implicit-explicit methods to nonlinear parabolic equations, providing optimal error estimates. Recent advances include the construction of new methods to solve stiff delay differential equations and study the stability and convergence of implicit-explicit methods, such as those proposed by Zhang and Xiao [25,26] and others.

In 1968, Alterman [35] first proposed using finite difference methods to simulate seismic wavefields, marking the beginning of its widespread application. Kelly [36] and others used finite differences to solve the 2D wave equation in 1976. Virieux [37,38] extended this to heterogeneous media with a second-order staggered grid method in 1984. Dablain [39] introduced higher-order difference methods, improving numerical accuracy. Moczo [41] and others in 2000 studied stability and grid discretization, identifying conditions for P and S wave stability. More recently, Kristek [42] introduced a method for discontinuous grid problems, significantly improving accuracy in 3D simulations.

This study innovatively divides the population into seven groups and establishes a reaction-diffusion model for the coronavirus, considering both time and space. The model is discretized in time using the Runge-Kutta and Euler implicit-explicit methods, and in space using the finite difference and finite element methods. Stability and convergence analyses of the discrete formats are conducted. The classic SIR model is extended by considering susceptible, isolated, exposed, asymptomatic infected, symptomatic infected, hospitalized, and recovered groups. Stability conditions and convergence results are provided, with model validation through numerical simulations.

## 2. Classic Infectious Disease SIR model

The SIR model divides the entire population into three parts, Where $S(t)$ It stands for susceptible population, which refers to people who are not sick but are easily infected; $I(t)$ The term "Infected population" refers to people who have been infected by the virus and can spread the virus to susceptible people; $R(t)$ refers to the population that has been rcovered from isolation or infection and has immunity. The specific model is as follows:

$$\frac{dS(t)}{dt} = -\beta S(t)I(t),$$
$$\frac{dI(t)}{dt} = \beta S(t)I(t) - \gamma I(t),$$
$$\frac{dR(t)}{dt} = \gamma I(t).$$

### 2.1. Novel Coronavirus Diffusion Model

Based on the transmission characteristics of the novel coronavirus and the theory of the classic SIR infectious disease model, a reaction-diffusion epidemic model for the novel coronavirus is established. The model is based on the following assumptions:

1. The parameters related to the model in the equations are all non-negative.
2. The model is a short-term model, and does not consider natural birth and death.
3. Patients who have recovered will not be reinfected in the short term.
4. Asymptomatic infections will not result in patient death.

let $x \in \Omega \subset R^d (d=1,2), t \in [0,T]$, the specific model is as follows：

$$
\begin{cases}
\dfrac{\partial Q(t,x)}{\partial t} = cS(t,x) - \delta Q(t,x), \\
\dfrac{\partial E(t,x)}{\partial t} = \theta S(t,x)(I(t,x) + bA(t,x)) - \varepsilon E(t,x) + \nabla \bullet (Nv_E \nabla E(t,x)), \\
\dfrac{\partial S(t,x)}{\partial t} = -\theta S(t,x)(I(t,x) + bA(t,x)) - cS(t,x) + \delta Q(t,x) + \nabla \bullet (Nv_S \nabla S(t,x)), \\
\dfrac{\partial I(t,x)}{\partial t} = \varepsilon f E(t,x) - jI(t,x) - lI(t,x) - h_1 I(t,x) + \nabla \bullet (Nv_I \nabla I(t,x)), \\
\dfrac{\partial D(t,x)}{\partial t} = gA(t,x) + h_1 I(t,x) - mD(t,x) - \mu D(t,x), \\
\dfrac{\partial R(t,x)}{\partial t} = \beta A(t,x) + jI(t,x) + \mu D(t,x),
\end{cases}
\quad (2.1)
$$

Among them，$N = S(x,t) + Q(x,t) + E(x,t) + A(x,t) + I(x,t) + D(x,t) + R(x,t)$ represent total population $\dfrac{\partial}{\partial t}$, represent differentiation in time direction, $\nabla$ represent gradient in space, parameter $\theta, b, c, f, g, \beta, \delta, \varepsilon, h_1, l, \mu, v_S, v_E, v_A, v_I$ The value range is $[0,1]$.

The dynamics of virus spread among seven different groups, namely susceptible people, isolated people, exposed people, asymptomatic infected people, symptomatic infected people, confirmed people in hospital and recovered people, are shown in Figure 2.1.

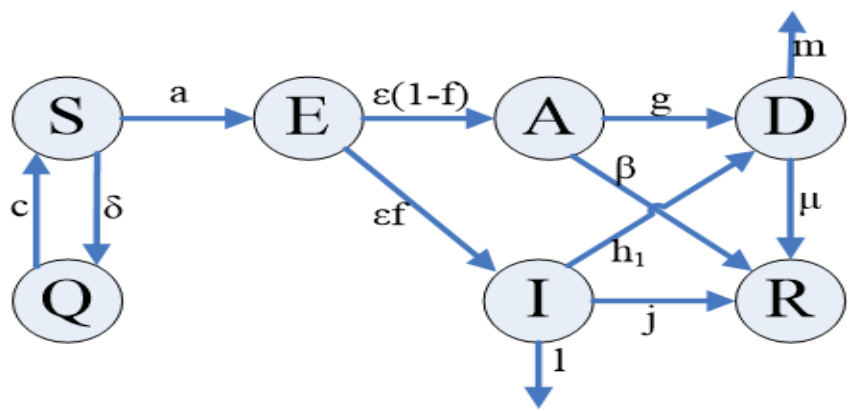

Figure 2.1 Dynamic diagram of infection among different populations

**Figure 2.1** shows the infection dynamics between the 7 different groups in the model. The exposed group $E(t,x)$ refers to individuals living in an environment with a low level of the virus, carrying a small amount of the virus but not being infectious. The proportion of this group that transitions into asymptomatic carriers is $\varepsilon(1 - f)$, and the proportion that transitions into symptomatic carriers is $\varepsilon f$. The isolated group $Q(t,x)$ refers to individuals who are in home quarantine or are isolated in designated locations as per national regulations. This group does not come into contact with the infected population and cannot be infected. The confirmed hospitalized treatment group $D(t,x)$ represents the number of confirmed cases receiving isolated treatment in government-designated hospitals at time **t**. A portion of individuals from the confirmed hospitalized treatment group and symptomatic carriers will die due to the infection, thereby being removed from the system.

## 2.2. COVID-19 Virus Diffusion Model: Runge-Kutta Finite Difference Discretization Format

This section discretizes model (2.1) in time using the Runge-Kutta method and in space using the finite difference method with $\Omega = [0, L]$. The grid is subdivided in the two-dimensional time-space region with the step size along the x-axis being $h = \frac{L}{N}$, and the step size along the t-axis being $\Delta t = TK$ with step size $x_i = ih, t_n = n\Delta t, i = 0,1,2,\ldots,N, n = 0,1,2,\ldots,K,$. The grid points are denoted as $(x_i, t_n)$.

First, the finite difference method is applied to discretize the space, resulting in the following semi-discrete format:

$$\frac{\partial S(t,x_i)}{\partial t} = -\theta S(t,x_i)\left(I(t,x_i) + bA(t,x_i)\right) - cS(t,x_i) + \delta Q(t,x_i)$$
$$+ Nv_S \frac{S(t,x_{i+1}) - 2S(t,x_i) + S(t,x_{i-1})}{h^2},$$

$$\frac{\partial Q(t,x_i)}{\partial t} = cS(t,x_i) - \delta Q(t,x_i),$$

$$\frac{\partial E(t,x_i)}{\partial t} = \theta S(t,x_i)\left(I(t,x_i) + bA(t,x_i)\right) - \varepsilon E(t,x_i)$$
$$+ Nv_E \frac{E(t,x_{i+1}) + E(t,x_{i-1}) - 2E(t,x_i)}{h^2},$$

$$\frac{\partial A(t,x_i)}{\partial t} = \varepsilon(1-f)E(t,x_i) - gA(t,x_i) - \beta A(t,x_i)$$
$$+ Nv_A \frac{A(t,x_{i+1}) + A(t,x_{i-1}) - 2A(t,x_i)}{h^2},$$

$$\frac{\partial I(t,x_i)}{\partial t} = \varepsilon f E(t,x_i) - jI(t,x_i) - lI(t,x_i) - h_1 I(t,x_i)$$
$$+ Nv_I \frac{I(t,x_{i+1}) + I(t,x_{i-1}) - 2I(t,x_i)}{h^2},$$

$$\frac{\partial D(t,x_i)}{\partial t} = gA(t,x_i) + h_1 I(t,x_i) - mD(t,x_i) - \mu D(t,x_i).$$

$$\frac{\partial R(t,x_i)}{\partial t} = \beta A(t,x_i) + jI(t,x_i) + \mu D(t,x_i).$$

Then, the model is time discretized using the Runge-Kutta method, and the fourth-order Runge-Kutta fully discrete format is as follows:

$$\begin{cases} S_{n+1,i} = S_{n,i} + \dfrac{\Delta t}{6}(k_{1i}^s + 2k_{2i}^s + 2k_{3i}^s + k_{4i}^s), \\ Q_{n+1,i} = Q_{n,i} + \dfrac{\Delta t}{6}(k_{1i}^q + 2k_{2i}^q + 2k_{3i}^q + k_{4i}^q), \\ E_{n+1,i} = E_{n,i} + \dfrac{\Delta t}{6}(k_{1i}^e + 2k_{2i}^e + 2k_{3i}^e + k_{4i}^e), \\ A_{n+1,i} = A_{n,i} + \dfrac{\Delta t}{6}(k_{1i}^a + 2k_{2i}^a + 2k_{3i}^a + k_{4i}^a), \\ I_{n+1,i} = I_{n,i} + \dfrac{\Delta t}{6}(k_{1i}^i + 2k_{2i}^i + 2k_{3i}^i + k_{4i}^i), \\ D_{n+1,i} = D_{n,i} + \dfrac{\Delta t}{6}(k_{1i}^d + 2k_{2i}^d + 2k_{3i}^d + k_{4i}^d), \\ R_{n+1,i} = R_{n,i} + \dfrac{\Delta t}{6}(k_{1i}^d + 2k_{2i}^d + 2k_{3i}^d + k_{4i}^d), \end{cases} \quad (2.2)$$

The first-order coefficients are as follows:

$$k_{1i}^s = -\theta S_{n,i}(I_{n,i} + bA_{n,i}) - cS_{n,i} + \delta Q_{n,i} + N v_S \frac{S_{n+1,i} - 2S_{n,i} + S_{n-1,i}}{h^2},$$

$$k_{1i}^q = cS_{n,i} - \delta Q_{n,i},$$

$$k_{1i}^e = \theta S_{n,i}(I_{n,i} + bA_{n,i}) - \varepsilon E_{n,i} + N v_E \frac{E_{n+1,i} - 2E_{n,i} + E_{n-1,i}}{h^2},$$

$$k_{1i}^a = \varepsilon(1-f)E_{n,i} - gA_{n,i} - \beta A_{n,i} + N v_A \frac{A_{n,i+1} - 2A_{n,i} + A_{n,i-1}}{h^2},$$

$$k_{1i}^i = \varepsilon f E_{n,i} - jI_{n,i} - lI_{n,i} - h_1 I_{n,i} + N v_I \frac{I_{n,i+1} - 2I_{n,i} + I_{n,i-1}}{h^2},$$

$$k_{1i}^d = gA_{n,i} + h_1 I_{n,i} - mD_{n,i} - \mu D_{n,i},$$

$$k_{1i}^r = \beta A_{n,i} + jI_{n,i} + \mu D_{n,i}.$$

The second-order coefficients are as follows:

$$k_{2i}^s = -\theta(S_{n,i} + \frac{\Delta t}{2}k_{1i}^s)\left(I_{n,i} + \frac{\Delta t}{2}k_{1i}^i + b\left(A_{n,i} + \frac{\Delta t}{2}k_{1i}^a\right)\right) - c(S_{n,i} + \frac{\Delta t}{2}k_{1i}^s)$$
$$+ \delta(Q_{n,i} + \frac{\Delta t}{2}k_{1i}^q) + N v_S \frac{S_{n+1,i} + \frac{\Delta t}{2}k_{1i+1}^s + S_{n-1,i} + \frac{\Delta t}{2}k_{1i-1}^s - 2S_{n,i} - \Delta t k_{1i}^s}{h^2}.$$

$$k_{2i}^q = c(S_{n,i} + \frac{\Delta t}{2}k_{1i}^s) - \delta(Q_{n,i} + \frac{\Delta t}{2}k_{1i}^q),$$

$$k_{2i}^e = \theta(S_{n,i} + \frac{\Delta t}{2}k_{1i}^s)\left(I_{n,i} + \frac{\Delta t}{2}k_{1i}^i + b\left(A_{n,i} + \frac{\Delta t}{2}k_{1i}^a\right)\right) - \varepsilon(E_{n,i} + \frac{\Delta t}{2}k_{1i}^e)$$

$$+ Nv_E \frac{E_{n+1,i} + \frac{\Delta t}{2}k^e_{1i+1} + E_{n-1,i} + \frac{\Delta t}{2}k^e_{1i-1} - 2E_{n,i} - \Delta t k^e_{1i}}{h^2},$$

$$k^a_{2i} = \varepsilon(1-f)(E_{n,i} + \frac{\Delta t}{2}k^e_{1i}) - g(A_{n,i} + \frac{\Delta t}{2}k^a_{1i}) - \beta(A_{n,i} + \frac{\Delta t}{2}k^a_{1i})$$

$$+ Nv_A \frac{A_{n,i+1} + \frac{\Delta t}{2}k^a_{1i+1} + A_{n,i-1} + \frac{\Delta t}{2}k^a_{1i-1} - 2A_{n,i} - \Delta t k^a_{1i}}{h^2},$$

$$k^i_{2i} = \varepsilon f(E_{n,i} + \frac{\Delta t}{2}k^e_{1i}) - j(I_{n,i} + \frac{\Delta t}{2}k^i_{1i}) - l(I_{n,i} + \frac{\Delta t}{2}k^i_{1i}) - h_1(I_{n,i} + \frac{\Delta t}{2}k^i_{1i})$$

$$+ Nv_I \frac{I_{n,i+1} + \frac{\Delta t}{2}k^i_{1i+1} + I_{n,i-1} + \frac{\Delta t}{2}k^i_{1i-1} - 2I_{n,i} - \Delta t k^i_{1i}}{h^2},$$

$$k^d_{2i} = g(A_{n,i} + \frac{\Delta t}{2}k^a_{1i}) + h_1(I_{n,i} + \frac{\Delta t}{2}k^i_{1i}) - m(D_{n,i} + \frac{\Delta t}{2}k^d_{1i}) - \mu(D_{n,i} + \frac{\Delta t}{2}k^d_{1i}),$$

$$k^r_{2i} = \beta(A_{n,i} + \frac{\Delta t}{2}k^a_{1i}) + j(I_{n,i} + \frac{\Delta t}{2}k^i_{1i}) + \mu(D_{n,i} + \frac{\Delta t}{2}k^d_{1i}).$$

The third-order coefficients are as follows:

$$k^s_{3i} = -\theta(S_{n,i} + \frac{\Delta t}{2}k^s_{2i})\left(I_{n,i} + \frac{\Delta t}{2}k^i_{2i} + b\left(A_{n,i} + \frac{\Delta t}{2}k^a_{2i}\right)\right) - c(S_{n,i} + \frac{\Delta t}{2}k^s_{2i}) + \delta(Q_{n,i} + \frac{\Delta t}{2}k^q_{2i})$$

$$+ Nv_S \frac{S_{n+1,i} + \frac{\Delta t}{2}k^s_{2i+1} + S_{n-1,i} + \frac{\Delta t}{2}k^s_{2i-1} - 2S_{n,i} - \Delta t k^s_{2i}}{h^2},$$

$$k^q_{3i} = c(S_{n,i} + \frac{\Delta t}{2}k^s_{2i}) - \delta(Q_{n,i} + \frac{\Delta t}{2}k^q_{2i}),$$

$$k^e_{3i} = \theta(S_{n,i} + \frac{\Delta t}{2}k^s_{2i})\left(I_{n,i} + \frac{\Delta t}{2}k^i_{2i} + b\left(A_{n,i} + \frac{\Delta t}{2}k^a_{2i}\right)\right) - \varepsilon(E_{n,i} + \frac{\Delta t}{2}k^e_{2i})$$

$$+ Nv_E \frac{E_{n+1,i} + \frac{\Delta t}{2}k^e_{2i+1} + E_{n-1,i} + \frac{\Delta t}{2}k^e_{2i-1} - 2E_{n,i} - \Delta t k^e_{2i}}{h^2},$$

$$k^a_{3i} = \varepsilon(1-f)(E_{n,i} + \frac{\Delta t}{2}k^e_{2i}) - g(A_{n,i} + \frac{\Delta t}{2}k^a_{2i}) - \beta(A_{n,i} + \frac{\Delta t}{2}k^a_{2i})$$

$$+ Nv_A \frac{A_{n,i+1} + \frac{\Delta t}{2}k^a_{2i+1} + A_{n,i-1} + \frac{\Delta t}{2}k^a_{2i-1} - 2A_{n,i} - \Delta t k^a_{2i}}{h^2},$$

$$k^i_{3i} = \varepsilon f(E_{n,i} + \frac{\Delta t}{2}k^e_{2i}) - j(I_{n,i} + \frac{\Delta t}{2}k^i_{2i}) - l(I_{n,i} + \frac{\Delta t}{2}k^i_{2i}) - h_1(I_{n,i} + \frac{\Delta t}{2}k^i_{2i})$$

$$+ Nv_I \frac{I_{n,i+1} + \frac{\Delta t}{2}k^i_{2i+1} + I_{n,i-1} + \frac{\Delta t}{2}k^i_{2i-1} - 2I_{n,i} - \Delta t k^i_{2i}}{h^2},$$

$$k_{3i}^d = g(A_{n,i} + \frac{\Delta t}{2} k_{2i}^a) + h_1(I_{n,i} + \frac{\Delta t}{2} k_{2i}^i) - m(D_{n,i} + \frac{\Delta t}{2} k_{2i}^d) - \mu(D_{n,i} + \frac{\Delta t}{2} k_{2i}^d),$$

$$k_{3i}^r = \beta(A_{n,i} + \frac{\Delta t}{2} k_{2i}^a) + j(I_{n,i} + \frac{\Delta t}{2} k_{2i}^i) + \mu(D_{n,i} + \frac{\Delta t}{2} k_{2i}^d).$$

The forth-order coefficients are as follows:

$$k_{4i}^s = -\theta(S_{n,i} + \Delta t k_{3i}^s)(I_{n,i} + \Delta t k_{3i}^i + b(A_{n,i} + \Delta t k_{3i}^a)) - c(S_{n,i} + \Delta t k_{3i}^s) + \delta(Q_{n,i} + \Delta t k_{3i}^q)$$
$$+ N v_S \frac{S_{n+1,i} + \Delta t k_{3i+1}^s + S_{n-1,i} + \Delta t k_{3i-1}^s - 2S_{n,i} - 2\Delta t k_{3i}^s}{h^2},$$

$$k_{3i}^q = c(S_{n,i} + \frac{\Delta t}{2} k_{2i}^s) - \delta(Q_{n,i} + \frac{\Delta t}{2} k_{2i}^q),$$

$$k_{4i}^e = \theta(S_{n,i} + \Delta t k_{3i}^s)(I_{n,i} + \Delta t k_{3i}^i + b(A_{n,i} + \Delta t k_{3i}^a)) - \varepsilon(E_{n,i} + \Delta t k_{3i}^e)$$
$$+ N v_E \frac{E_{n+1,i} + \Delta t k_{3i+1}^e + E_{n-1,i} + \Delta t k_{3i-1}^e - 2E_{n,i} - 2\Delta t k_{3i}^e}{h^2},$$

$$k_{4i}^a = \varepsilon(1-f)(E_{n,i} + \Delta t k_{3i}^e) - g(A_{n,i} + \Delta t k_{3i}^a) - \beta(A_{n,i} + \Delta t k_{3i}^a)$$
$$+ N v_A \frac{A_{n,i+1} + \Delta t k_{3i+1}^a + A_{n,i-1} + \Delta t k_{3i-1}^a - 2A_{n,i} - 2\Delta t k_{3i}^a}{h^2},$$

$$k_{4i}^i = \varepsilon f(E_{n,i} + \Delta t k_{3i}^e) - j(I_{n,i} + \Delta t k_{3i}^i) - l(I_{n,i} + \Delta t k_{3i}^i) - h_1(I_{n,i} + \Delta t k_{3i}^i)$$
$$+ N v_I \frac{I_{n,i+1} + \Delta t k_{3i+1}^i + I_{n,i-1} + \Delta t k_{3i-1}^i - 2I_{n,i} - 2\Delta t k_{3i}^i}{h^2},$$

$$k_{4i}^d = g(A_{n,i} + \Delta t k_{3i}^a) + h_1(I_{n,i} + \Delta t k_{3i}^i) - m(D_{n,i} + \Delta t k_{3i}^d) - \mu(D_{n,i} + \Delta t k_{3i}^d),$$

$$k_{4i}^r = \beta(A_{n,i} + \Delta t k_{3i}^a) + j(I_{n,i} + \Delta t k_{3i}^i) + \mu(D_{n,i} + \Delta t k_{3i}^d).$$

## 2.3 Stability Analysis of the Runge-Kutta Finite Difference Format for the COVID-19 Virus Diffusion Model

This section presents the stability analysis of the Runge-Kutta finite difference format for the COVID-19 virus diffusion model. First, the following lemma is introduced:

Holder inequality applies to any $f \in L^p(\Omega)$ and $g \in L^q(\Omega)$ gives,

$$\|fg\|_{L^1(\Omega)} \leq \|f\|_{L^p(\Omega)} \|g\|_{L^q(\Omega)}.$$

with $1 \leq p, q \leq \infty$, satisfies $\frac{1}{p} + \frac{1}{q} = 1$,

**Schwarz** inequality, when $p = q = 2$, **Holder** inequality is:

$$\|fg\|_{L^1(\Omega)} \leq \|f\|_{L^2(\Omega)} \|g\|_{L^2(\Omega)}.$$

**Gronwall** lemma assumes that $C_0, \Delta t$ is a positive integer, and $a_k, b_k, c_k, d_k$ is non negative sequence that satisfy the following condition:

$$a_n + \Delta t \sum_{k=0}^{n} b_k \leq \Delta t \sum_{k=0}^{n-1} d_k a_k + \Delta t \sum_{k=0}^{n-1} c_k + C_0, \quad \forall n \geq 1,$$

That is,

$$a_n + \Delta t \sum_{k=0}^{n} b_k \leq \left( \Delta t \sum_{k=0}^{n-1} c_k + C_0 \right) \exp\left( \Delta t \sum_{k=0}^{n-1} d_k \right), \quad \forall n \geq 1$$

Let $y = [S, Q, E, A, I, D, R]^T$ define the function $f(t, y)$

$$f(t,y) = \begin{bmatrix} -\theta S(t_n, x_i)(I(t_n, x_i) + bA(t_n, x_i)) - cS(t_n, x_i) + \delta Q(t_n, x_i) + Nv_S \frac{S(t_n, x_{i+1}) + S(t_n, x_{i-1}) - 2S(t_n, x_i)}{h^2} \\ cS(t_n, x_i) - \delta Q(t_n, x_i) \\ \theta S(t_n, x_i)(I(t_n, x_i) + bA(t_n, x_i)) - \varepsilon E(t_n, x_i) + Nv_E \frac{E(t_n, x_{i+1}) + E(t_n, x_{i-1}) - 2E(t_n, x_i)}{h^2} \\ \varepsilon(1-f)E(t_n, x_i) - gA(t_n, x_i) - \beta A(t_n, x_i) + Nv_A \frac{A(t_n, x_{i+1}) + A(t_n, x_{i-1}) - 2A(t_n, x_i)}{h^2} \\ \varepsilon f E(t_n, x_i) - jI(t_n, x_i) - lI(t_n, x_i) - h_1 I(t_n, x_i) + Nv_I \frac{I(t_n, x_{i+1}) + I(t_n, x_{i-1}) - 2I(t_n, x_i)}{h^2} \\ gA(t_n, x_i) + h_1 I(t_n, x_i) - mD(t_n, x_i) - \mu D(t_n, x_i) \\ \beta A(t_n, x_i) + jI(t_n, x_i) + \mu D(t_n, x_i) \end{bmatrix},$$

By [43], we have the following lemma

Lemma 2.1 Let $\langle \bullet \ \bullet \rangle$ denote the inner product on n-dimensional space, $\|\bullet\|$ be the norm on $n$ dimensional space $R^n$, for $\forall t \geq 0, y_1, y_2 \in R^n$, we have
$$\langle y_1 - y_2, f(t, y_1) - f(t, y_2) \rangle \leq \alpha \| y_1 - y_2 \|^2$$
Where $\alpha$ is any constant number.

The following presents the stability analysis of the **Runge-Kutta finite difference discretization format (2.2)**.

**Theorem 2.1**: Let the solution of the Runge-Kutta finite difference discretization format (2.2) for model (2.1) be denoted as $y_{n,i}$. When the time step size $\Delta t$ satisfies $\Delta t \leq 1/\alpha$, the diffusion model's Runge-Kutta finite difference discretization format is stable, with the following estimate:
$$\| y_{n,i} \|^2 \leq \| y_{0,i} \|^2,$$

where $\alpha$ is nonzero positive number.

Proof: Model (2.1) can be expressed as
$$y' = f(t, y).$$

Implies,

$$Y_{j,i}^{(n)} = y_{n,i} + \Delta t \sum_{k=1}^{4} a_{jk} f(t_n + c_k \Delta t, Y_{k,i}^{(n)}), \quad j=1,2,3,4,$$

$$y_{n+1,i} = y_{n,i} + \Delta t \sum_{k=1}^{4} b_k f(t_n + c_k \Delta t, Y_{k,i}^{(n)}),$$

$$Y_{j,i}^{(n)} = y(t_n + c_j \Delta t, xi),$$

$$Q_{k,i}^{(n)} = f(t_n + c_k \Delta t, Y_{k,i}^{(n)}),$$

and

$$y_{n,i} = y(t_n, x_i), t_n = n\Delta t,$$

$$a_{jk} = \begin{pmatrix} 0 & 0 & 0 & 0 \\ 1/2 & 0 & 0 & 0 \\ 0 & 1/2 & 0 & 0 \\ 0 & 0 & 0 & 0 \end{pmatrix},$$

$$b_k = (1/6 \quad 1/3 \quad 1/3 \quad 1/6),$$

$$c_k = (0 \quad 1/2 \quad 1/2 \quad 1), j=1,2,3,4; k=1,2,3,4.$$

we get,

$$\|y_{n+1,i}\|^2 - \|y_{n,i}\|^2 - 2\sum_{k=1}^{4} d_k \langle Y_{k,i}^{(n)}, \Delta t Q_{k,i}^{(n)} - Y_{k,i}^{(n)} \rangle = -\sum_{j=1}^{5}\sum_{k=1}^{5} M_{jk} \langle \xi_j, \xi_k \rangle, \quad (2.3)$$

There exist,

$$\xi_1 = y_{n,i}, \xi_j = \Delta t Q j - 1, in, j=2,3,4,5,$$

$$M_{jk} = \begin{bmatrix} -2e_4^T D_4 e_4 & e_4^T D_4 - b_k^T - 2e_4^T D_4 a_{jk} \\ D_4 e_4 - b_k - 2a_{jk}^T D_4 e_4 & D_4 a_{jk} + a_{jk}^T D_4 - b_k b_k^T - 2a_{jk}^T D_4 a_{jk} \end{bmatrix},$$

$$D_4 = \begin{bmatrix} d_1 & 0 & 0 & 0 \\ 0 & d_2 & 0 & 0 \\ 0 & 0 & d_3 & 0 \\ 0 & 0 & 0 & d_4 \end{bmatrix}.$$

where $D_4$ is non negative diagonal matrix.

$$e_4 = [1,1,1,1]^T$$

Because,

$$-\sum_{j=1}^{5}\sum_{k=1}^{5} M_{jk} \langle y_{n,i}, \Delta t Q_{j-1,i}^{(n)} \rangle \leq 0,$$

Equation (2.3) can be expressed as:

$$\|y_{n+1,i}\|^2 - \|y_{n,i}\|^2 - 2\sum_{k=1}^{4} d_k \langle Y_{k,i}^{(n)}, \Delta t Q_{k,i}^{(n)} - Y_{k,i}^{(n)} \rangle \leq 0.$$

Therefore,

$$\|y_{n+1,i}\|^2 \leq \|y_{n,i}\|^2 + 2\sum_{k=1}^{4} d_k \langle Y_{k,i}^{(n)}, \Delta t Q_{k,i}^{(n)} - Y_{k,i}^{(n)} \rangle. \qquad (2.4)$$

Because,

$$\exists Y_{0,i}^{(n)} = 0, st.\ f(t, Y_{0,i}^{(n)}) = 0,$$

Gives:

$$2 \langle Y_{k,i}^{(n)}, \Delta t Q_{k,i}^{(n)} \rangle = 2 \langle Y_{k,i}^{(n)}, f(t_n + c_k \Delta t, Y_{k,i}^{(n)}) \rangle = 2 \langle Y_{k,i}^{(n)} - Y_{0,i}^{(n)}, f(t_n + c_k \Delta t, Y_{k,i}^{(n)}) \rangle$$
$$= 2 \langle Y_{k,i}^{(n)} - Y_{0,i}^{(n)}, f(t_n + c_k \Delta t, Y_{k,i}^{(n)}) - f(t_n + c_k \Delta t, Y_{0,i}^{(n)}) \rangle.$$

From lemma 2.1, we can say that:

$$2 \langle Y_{k,i}^{(n)}, \Delta t Q_{k,i}^{(n)} \rangle \leq 2\Delta t \alpha \|Y_{k,i}^{(n)} - Y_{0,i}^{(n)}\|^2 = 2\Delta t \alpha \|Y_{k,i}^{(n)}\|^2 \qquad (2.5)$$

By substituting (2.5) into (2.4), we get:

$$\|y_{n+1,i}\|^2 \leq \|y_{n,i}\|^2 + 2\sum_{k=1}^{4} d_k \langle Y_{k,i}^{(n)}, \Delta t Q_{k,i}^{(n)} - Y_{k,i}^{(n)} \rangle \qquad (2.6)$$
$$\leq \|y_{n,i}\|^2 + (2\Delta t \alpha - 2)\sum_{k=1}^{4} d_k \|Y_{k,i}^{(n)}\|$$

Because,

$$\Delta t \leq 1/\alpha,$$

We get,

$$2\Delta t \alpha - 2 \leq 0.$$

Equation (2.6) can be transformed into:

$$\|y_{n+1,i}\|^2 \leq \|y_{n,i}\|^2$$

Then

$$\|y_{n+1,i}\|^2 \leq \|y_{0,i}\|^2$$

This completes the proof.

## 2.4 Convergence analysis of the Runge-Kutta finite difference scheme of the new coronavirus diffusion model

This section will conduct a convergence analysis on the Runge-Kutta finite difference scheme of the new coronavirus diffusion model.

**Theorem 2.2** let $\tilde{y}(t_n, x_i)$ is true solution of model (2.1), $y_{n,i}$ is the solution od Runge-Kutta finite difference discretization scheme(2.2) of model(2.1) and $e_i^n = \tilde{y}(t_n, x_i) - y_{n,i}$, $0 \leq i \leq J$, $0 \leq n \leq N$, under the condition of Theorem 2.1 if the following estimation equationis satisfied:

Then the Runge-Kutta finite difference discrete scheme solution of model (2.1) is convergent.
$$\|e_i^n\| \leq C(h^2 + (\Delta t)^4),$$
Proof: The Runge-Kutta finite difference discrete scheme (2.2) of model (2.1) can be expressed as

$$y_{n+1,i} = y_{n,i} + \Delta t \sum_{k=1}^{4} b_k f(t_n + c_k \Delta t, Y_{k,i}^{(n)}),$$

let

$$\frac{y_{n+1,i} - y_{n,i}}{\Delta t} = \sum_{k=1}^{4} b_k f(t_n + c_k \Delta t, Y_{k,i}^{(n)}),$$

and

$$Y_{k,i}^{(n)} = y(t_n + c_k \Delta t, x_i).$$

Entering the true solution and the discrete format solution can be obtained

$$\frac{\tilde{y}(t_{n+1}, x_i) - \tilde{y}(t_n, x_i)}{\Delta t} = \sum_{k=1}^{4} b_k f(t_n + c_k \Delta t, \tilde{y}(t_n + c_j \Delta t, x_i)) + R_i^n, \quad (2.7)$$

$$\frac{y_{n+1,i} - y_{n,i}}{\Delta t} = \sum_{k=1}^{4} b_k f(t_n + c_k \Delta t, y(t_n + c_j \Delta t, x_i)). \quad (2.8)$$

where $R_i^n$ is truncation error $R_i^n = O(h^2 + \Delta t^4)$.

by substracting (2.8) from (2.8), we can obtain

$$\frac{e_i^{n+1} - e_i^n}{\Delta t} = \sum_{k=1}^{4} b_k \{f(t_n + c_k \Delta t, \tilde{y}(t_n + c_j \Delta t, x_i)) - f(t_n + c_k \Delta t, y(t_n + c_j \Delta t, x_i))\} + R_i^n. \quad (2.9)$$

hence

$$F = \sum_{k=1}^{4} b_k \{f(t_n + c_k \Delta t, \tilde{y}(t_n + c_j \Delta t, x_i)) - f(t_n + c_k \Delta t, y(t_n + c_j \Delta t, x_i))\}.$$

Where as

By inner product of (2.9) and $e_i^{n+1} + e_i^n$, we can obtain:

$$\frac{\|e_i^{n+1}\|^2 - \|e_i^n\|^2}{\Delta t} = \langle F + R_i^n, e_i^{n+1} + e_i^n \rangle,$$

and

$$\begin{aligned}\langle F + R_i^n, e_i^{n+1} + e_i^n \rangle &= \langle F, e_i^{n+1} + e_i^n \rangle + \langle R_i^n, e_i^{n+1} + e_i^n \rangle \\ &\leq C\|F\|^2 + C(\|e_i^{n+1}\|^2 + \|e_i^n\|^2) + C(h^2 + (\Delta t)^4)^2 \\ &\leq C(h^2 + (\Delta t)^4)^2 + C(\|e_i^{n+1}\|^2 + \|e_i^n\|^2).\end{aligned}$$

We get,

$$\|e_i^{n+1}\|^2 - \|e_i^n\|^2 \leq \Delta t C(\|e_i^{n+1}\|^2 + \|e_i^n\|^2 + C(h^2 + (\Delta t)^4)^2).$$

By Gronwall lemma, we get

$$\|e_i^n\| \leq C(h^2 + (\Delta t)^4).$$

## 2.5 Numerical Examples

In this section, we select Jiangsu Province as the research object, use the finite difference method for spatial discretization, and the Runge-Kutta method for time discretization, and divide the overall administrative area of Jiangsu Province into quadrilateral grid areas. Through numerical simulation calculations, the predicted data is compared with the actual data to verify the effectiveness of the model.

$$\begin{cases} \dfrac{\partial S(t,x)}{\partial t} = -\theta S(t,x)(I(t,x) + bA(t,x)) - cS(t,x) + \delta Q(t,x), \\ \dfrac{\partial Q(t,x)}{\partial t} = cS(t,x) - \delta Q(t,x), \\ \dfrac{\partial E(t,x)}{\partial t} = \theta S(t,x)(I(t,x) + bA(t,x)) - \varepsilon E(t,x), \\ \dfrac{\partial A(t,x)}{\partial t} = \varepsilon(1-f)E(t,x) - gA(t,x) - \beta A(t,x), \\ \dfrac{\partial I(t,x)}{\partial t} = \varepsilon f E(t,x) - jI(t,x) - lI(t,x) - h_1 I(t,x), \\ \dfrac{\partial D(t,x)}{\partial t} = gA(t,x) + h_1 I(t,x) - mD(t,x) - \mu D(t,x), \\ \dfrac{\partial R(t,x)}{\partial t} = \beta A(t,x) + jI(t,x) + \mu D(t,x). \end{cases}$$

If the spatial spread of susceptible groups, exposed groups, asymptomatic infected persons, and symptomatic infected persons is not considered, that is, when the diffusion coefficients $v_E, v_S, v_A$ and $v_I$ and are all are 0, the new coronavirus diffusion equation becomes an ordinary differential equation. The specific format is as follows

Taking Nanjing City, Jiangsu Province as an example, a numerical simulation was conducted. The numerical values of each group of people in the model correspond to Jiangsu Province's 2021 City population on July 21. The numerical simulation results are shown in Figure 2.2.

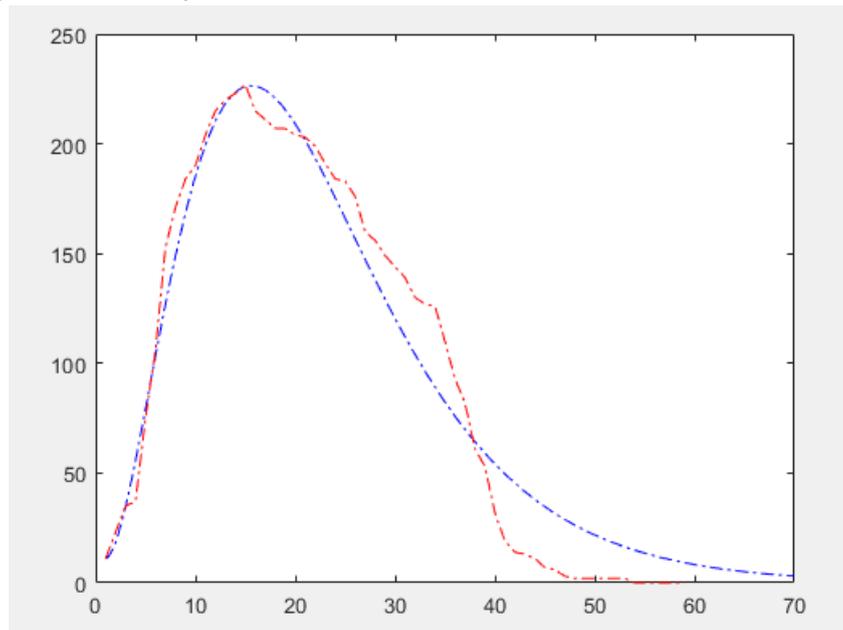

Figure 2.2 Comparison between actual data and Forecast data in Nanjing Area

**Figure 2.2** shows that the number of confirmed hospitalized cases in Nanjing reached its peak on the 15th day of the outbreak. Under strict epidemic prevention and control measures, the number of infections gradually decreased, and the number of confirmed hospitalized cases progressively approached zero. The numerical simulation results indicate that, without considering the spatial spread of the COVID-19 virus, the model's predicted results align well with the actual data from Nanjing's July epidemic, thereby validating the effectiveness of the model. When considering the spatial spread of susceptible, exposed, asymptomatic infected, and symptomatic infected individuals, Jiangsu Province is still taken as an example. At this point, the diffusion coefficients $v_S, v_E, v_A$ and $v_I$ satisfy:

The results are shown in figure 2.3

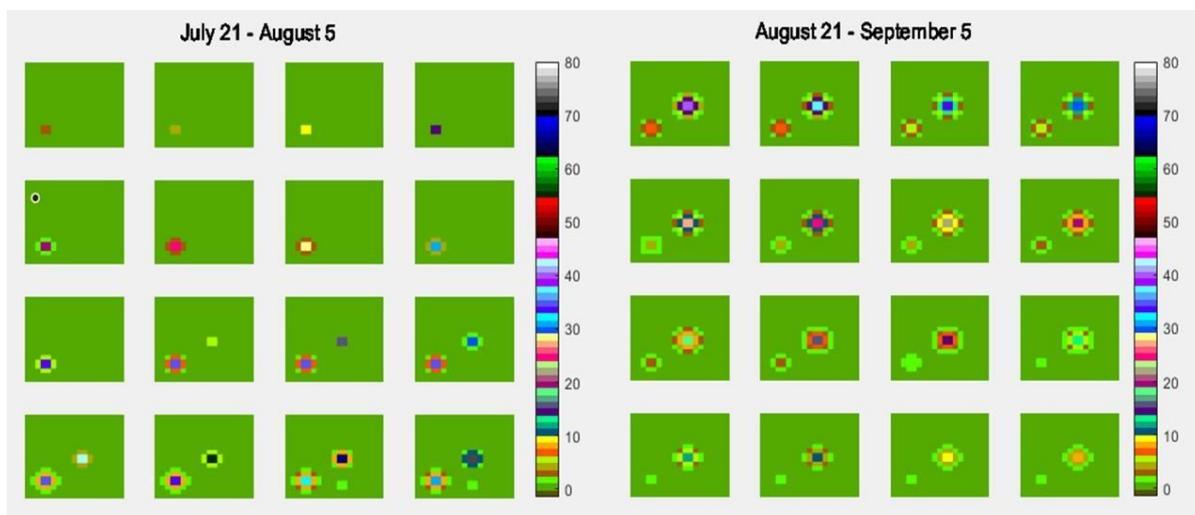

Figure 2.3 Changes in cases in various regions of Jiangsu Province.

**Figure 2.3** shows the spread of the COVID-19 pandemic in Jiangsu Province since July 21, 2021. On July 21, 2021, an outbreak of COVID-19 occurred in Nanjing. With the movement of people, the epidemic spread to cities such as Yangzhou and Huai'an. During this period, the epidemic spread rapidly, and Yangzhou became the most severely affected city in Jiangsu Province. From Figure 2.3, it can be seen that although the outbreaks in Nanjing and Yangzhou occurred relatively quickly, strict epidemic prevention policies helped contain the spread of the virus, preventing a massive surge in other cities within Jiangsu Province. The numerical simulation results are in good agreement with the actual situation, validating the effectiveness of the model.